\documentclass[a4paper]{amsart}

\usepackage{epic, eepic, amsfonts, latexsym, amssymb, graphicx,
multicol, mathrsfs, color, amscd, verbatim, paralist,
xspace, url, euscript, stmaryrd,  amsmath, enumitem,
bbold, multirow, tikz, mathtools}
\usepackage[all,pdf,cmtip]{xy}

\usepackage[colorlinks, linkcolor=blue, citecolor=magenta, urlcolor=cyan]{hyperref}

\usepackage{graphicx}

\usepackage{float}

\renewcommand\bar{\overline}

\def\RR{{\mathbb R}}
\def\CC{{\mathbb C}}

\newcommand{\suchthat}{\;\ifnum\currentgrouptype=16 \middle\fi|\;}

\def\Im{\mathop{\rm Im}\nolimits}
\def\Re{\mathop{\rm Re}\nolimits}

\def\h{\mathop{\mathfrak {hol}}\nolimits}

\def\su{\mathop{\mathfrak{su}}\nolimits}
\def\sl{\mathop{\mathfrak{sl}}\nolimits}

\def\qed{{\hfill $\Box$}}
\newtheorem{theorem}{THEOREM}[section]
\newtheorem{corollary}[theorem]{Corollary}

\newtheorem{conjecture}[theorem]{Conjecture}

\theoremstyle{definition}

\theoremstyle{remark}
\newtheorem{remark}[theorem]{Remark}

\makeatletter
\def\blfootnote{\xdef\@thefnmark{}\@footnotetext}
\makeatother

\begin{document}

\title[A short proof of the Dimension Conjecture]{A short proof of the Dimension Conjecture
\vspace{0.1cm}\\
for real hypersurfaces in $\CC^2$}
\author[Isaev]{Alexander Isaev}
\author[Kruglikov]{Boris Kruglikov}

\address[Isaev]{Mathematical Sciences Institute\\
Australian National University\\
Acton, ACT 2601, Australia}
\email{alexander.isaev@anu.edu.au}

\address[Kruglikov]{Department of Mathematics and Statistics\\
University of Troms\o{}\\
Troms\o{} 90-37, Norway}
\email{boris.kruglikov@uit.no}

\maketitle

\thispagestyle{empty}

\pagestyle{myheadings}

\begin{abstract} Recently, I.~Kossovskiy and R.~Shafikov have settled the so-called Dimension Conjecture, which characterizes spherical hypersurfaces in $\CC^2$ via the dimension of the algebra of infinitesimal automorphisms. In this note, we propose a short argument for obtaining their result.
\end{abstract}

\section{Introduction}\label{intro}
\setcounter{equation}{0}

Let $M$ be a 3-dimensional connected real-analytic CR-manifold of hypersurface type. We only consider $M$ locally, and therefore one can assume that $M$ is embedded in $\CC^2$ with the CR-structure induced by the complex structure of the ambient space. Recall that an infinitesimal CR-automorphism of $M$ is a smooth vector field on $M$ whose flow consists of CR-transformations. For $p\in M$, denote by $\h(M,p)$ the Lie algebra of real-analytic infinitesimal CR-automorphisms of $M$ defined in a neighborhood of $p$ on $M$, with the neighborhood {\it a priori}\, depending on the vector field. It is not hard to show that every element of $\h(M,p)$ is the real part of a holomorphic vector field defined on an open subset of $\CC^2$. 

If $M$ is Levi-flat, i.e., its Levi form identically vanishes, then $M$ is locally CR-equivalent to the direct product $\CC\times\RR\subset\CC^2$, hence in this case $\dim\h(M,p)=\infty$ for all $p\in M$. On the other hand, if $M$ is Levi nondegenerate at some point, then for every $p\in M$ one has $\dim\h(M,p)<\infty$ (see \cite[Theorem 11.5.1 and Corollary 12.5.5]{BER}). If, furthermore, $M$ is spherical at $p$, i.e., CR-equivalent to an open subset of the sphere $S^3\subset\CC^2$ in a neighborhood of $p$, then $\dim\h(M,p)=8$. Indeed, for every $q\in S^3$ the algebra $\h(S^3,q)$ consists of globally defined vector fields and is isomorphic to $\su_{2,1}$ (see \cite{P}, \cite{C} as well as \cite{CM}, \cite{Ta}, \cite[pp.~211--219]{Sa}, \cite{B2} for generalizations to higher CR-dimensions and CR-codimensions). Further, the reduction of 3-dimensional Levi nondegenerate CR-structures to absolute parallelisms obtained by \'E.~Cartan in \cite{C} implies that 8 is the maximal possible value of $\dim\h(M,p)$ provided the Levi form of $M$ at $p$ does not vanish. Moreover, as noted in \cite[part I, p.~34]{C}, results of A.~Tresse in \cite{Tr} (see also \cite{K}, \cite{KT}) yield that for such a point $p$ the condition $\dim\h(M,p)>3$ forces $M$ to be spherical at $p$.

This note concerns the following conjecture, in which Levi nondegeneracy is no longer assumed:

\begin{conjecture}\label{conj1} If $M$ is not Levi-flat, then for any $p\in M$ the condition
\begin{equation}
\dim\h(M,p)>5\label{thecondition}
\end{equation}
implies that $M$ is spherical at $p$.\end{conjecture} 

\noindent In the above form, the conjecture was formulated in article \cite{KS} where the authors called it the Dimension Conjecture and argued that it can be viewed as a variant of Poincar\'eÕs {\it probl\`eme local}. This statement is also a refined version, in the case $n=2$, of another conjecture, due to V.~Beloshapka, proposed in \cite{B3} for real hypersurfaces in $\CC^n$ with any $n\ge 2$, which so far has only been resolved for $n\le 3$ (see \cite{IZ}). 

In \cite{KS}, the authors proved:

\begin{theorem}\label{main} Conjecture {\rm \ref{conj1}} holds true.
\end{theorem}

\noindent The method of \cite{KS} is rather involved and based on considering second-order complex ODEs with meromorphic singularity. The aim of the present paper is to provide a short proof of Theorem \ref{main} by using standard facts on Lie algebras and their actions. Before proceeding, we state the following:

\begin{corollary}\label{cor}
The possible dimensions of $\h(M,p)$ are {\rm 0, 1, 2, 3, 4, 5, 8,} $\infty$, and each of these possibilities is realizable.
\end{corollary} 

{\bf Acknowledgements.} We thank B. Lamel for useful discussions and V. Beloshapka for communicating to us an example of a hypersurface in $\CC^2$ with\linebreak $\dim \h(M,p)=1$, for which he kindly performed {\tt Maple}-assisted computations. The first author is supported by the Australian Research Council.

\section{Proof of Theorem \ref{main} and Corollary \ref{cor}}\label{sect1}
\setcounter{equation}{0}

Suppose that $M$ is not Levi-flat. Then the set ${\mathcal S}$ of points of Levi nondegeneracy is dense in $M$. Fix $p\in M$ with $\dim\h(M,p)> 5$ and consider the algebra $\h(M,p)$. If $p\in {\mathcal S}$, then, as stated in the introduction, the sphericity of $M$ at $p$ follows from classical results in \cite{C}, \cite{Tr}.

Assume now that $p\not\in{\mathcal S}$. As $\dim\h(M,p)<\infty$, there exists a neighborhood $U$ of $p$ in $M$ where all vector fields in $\h(M,p)$ are defined. Therefore, for every\linebreak $p'\in U\cap{\mathcal S}$, the algebra $\h(M,p)$ is a subalgebra of $\h(M,p')$. Arguing as above, we then see that $M$ is spherical at $p'$. Hence, $\h(M,p)$ can be identified with a subalgebra of $\su_{2,1}$. It is not hard to show that $\su_{2,1}$ has no subalgebras of dimensions 6 and 7. This is a consequence, for instance, of the proof of Lemma 2.4 in \cite{EaI}, but for the reader's convenience we give a different argument here. Indeed, by \cite{M}, a maximal proper subalgebra of a semi-simple Lie algebra is either parabolic, or semi-simple or the stabilizer of a pseudo-torus. Therefore, all maximal subalgebras of $\su_{2,1}$ up to conjugation are described as follows: (i) one parabolic subalgebra, of dimension 5;\ (ii) one semi-simple subalgebra, namely $\mathfrak{so}_{2,1}$, of dimension 3; (iii) two pseudotoric subalgebras, namely $\mathfrak{u}_2$ and $\mathfrak{u}_{1,1}$, both of dimension 4. In particular, $\su_{2,1}$ has no subalgebras of dimension 6 and 7 as claimed.

Thus, we have $\h(M,p)=\su_{2,1}$. Consider the isotropy subalgebra $\h_0(M,p)\subset\h(M,p)$, which consists of all vector fields in $\h(M,p)$ vanishing at $p$. Clearly, $\dim\h_0(M,p)\ge 5$, and we obtain, again by the nonexistence of codimension one and two subalgebras in $\su(2,1)$, that either $\dim\h_0(M,p)=5$ or $\dim\h_0(M,p)=8$. In the former case, it follows that the orbit of $p$ under the corresponding local action of $SU(2,1)$ is open. Since $M$ is spherical at every point $p'\in U\cap{\mathcal S}$, we then see that $M$ is spherical at $p$ as required.

Suppose now that $\dim\h_0(M,p)=8$, i.e., $\h_0(M,p)=\su_{2,1}$. As shown in \cite{GS} (see pp.~113--115 therein), an action of a semisimple Lie algebra ${\mathfrak g}$ by real-analytic vector fields on a real-analytic manifold $X$ can be linearized near a fixed point $x$, i.e., there exist local coordinates in a neighborhood of $x$ on $X$ in which all vector fields arising from ${\mathfrak g}$ are linear. It then follows that $\su_{2,1}$ has a nontrivial real 3-dimensional representation. On the other hand, it is easy to see that no such representation exists. Indeed, assuming the contrary and complexifying, we obtain a complex 3-dimensional representation of $\sl_3(\CC)$. Up to isomorphism, this is the standard (defining) representation, hence the standard action of $\su_{2,1}$ on $\CC^3$ must have an invariant totally real 3-dimensional subspace, and it is straightforward to verify that no such subspace in fact exists. This contradiction completes the proof of the theorem.\qed        

\begin{remark}\label{8dim} The argument contained in the last paragraph of the above proof provides a short way of answering the question asked in the title of article \cite{B4}.
\end{remark}

Next, to prove Corollary \ref{cor}, we only need to observe that each of the integers 0, 1, 2, 3, 4, 5 is realizable as $\dim\h(M,p)$. The realizability of 0, 2, 3, 4, 5 follows from the examples given in \cite[p.~143]{B4}, \cite[Table 1]{KL}, \cite{St}, so it only remains to find an example with $\dim\h(M,p)=1$. Consider the hypersurface $\Gamma_1$ given in coordinates $z,w$ in $\CC^2$ by the equation\footnote{This example was communicated to us by V.~Beloshapka.}
$$
\Im w=|z|^2+(\Re z^2)|z|^2.
$$
By Theorem 3 of \cite{B1}, the stability group of $\Gamma_1$ at the origin consists only of the transformations $z\mapsto \pm z, w\mapsto w$, hence $\h_0(\Gamma_1,0)=0$. One can further show (e.g., by {\tt Maple}-assisted computations) that $\h(\Gamma_1,0)$ is spanned by the vector field $\partial/\partial w+\partial/\partial\bar{w}$. Another example is given by the hypersurface $\Gamma_2$ defined as
$$
\Im w=|z|^2+(\Re w) |z|^8.
$$
In this case, the stability group at the origin consists of all rotations in $z$ (see, e.g., \cite[p.~1159]{EzI}), and one can further show that every element of $\h(\Gamma_2,0)$ vanishes at the origin. Hence, $\h(\Gamma_2,0)$ is spanned by $iz\partial/\partial z-i\bar{z}\partial/\partial\bar{z}$. One can produce many more examples of this kind by considering hypersurfaces of the form $\Im w=f(|z|^2,\Re w)$, where $f$ is real-analytic and in general position.\qed

\begin{remark}\label{5dim}  As we noted in the proof of Theorem \ref{main}, $\su_{2,1}$ has only one, up to conjugation, 5-dimensional subalgebra (which is parabolic), and this is exactly the subalgebra that occurs in the examples with $\dim\h(M,p)=5$ given in \cite{B4}, \cite{KL}. In all these cases, one has $\h(M,p)=\h_0(M,p)$. Explicitly classifying the manifolds with $\dim\h(M,p)=5$ requires a much greater effort, and article \cite{KS} makes progress in this direction by showing that every such manifold has to be a \lq\lq sphere blowup\rq\rq\, (as defined above the statement of Theorem 3.10 therein). 
\end{remark}

\end{document}